\input amstex
 \documentstyle{amsppt}
\NoBlackBoxes
\magnification1200
\pagewidth{6.5 true in} 
\pageheight{9.25 true in}
\topmatter
\title
Small gaps between prime numbers: the work of Goldston-Pintz-Y\i ld\i r\i m
\endtitle  
\author K. Soundararajan 
\endauthor 
\rightheadtext{Small gaps between primes} 
\thanks 
The author is partially supported by the National Science Foundation.
\endthanks
\address 
Department of Mathematics, University of Michigan, Ann Arbor, MI 48109, USA
\endaddress
\email 
ksound\@ \!\!umich.edu
\endemail
\endtopmatter
\def\lam{\lambda}
\def\Yildirim{Y{\i}ld{\i}r{\i}m\ }
\document
\loadbold
 
\noindent{\bf Introduction.} In early 2005, Dan Goldston, J{\' a}nos Pintz, and Cem \Yildirim [12]
 made a spectacular breakthrough in the study of prime numbers.  
 Resolving a long-standing open problem, they proved that there are infinitely many primes for which 
 the gap to the next prime is as small as we want compared to the average 
 gap between consecutive primes.  
 Before their work, it was only known that there were infinitely many gaps which were 
 about a quarter the size of the average gap.  The 
 new result may be viewed as a step towards the famous twin prime conjecture 
 that there are infinitely many prime pairs $p$ and $p+2$; the gap here 
 being $2$, the smallest possible gap between primes \footnote{apart from the gap between $2$ 
 and $3$, of course!}.   Perhaps most excitingly, their work reveals a connection between 
 the distribution of primes in arithmetic progressions and small gaps between primes.
 Assuming certain (admittedly difficult) conjectures on the distribution of primes 
 in arithmetic progressions, they are able to prove the existence of infinitely many prime 
 pairs that differ by at most $16$.  The aim of this article is to explain some of the ideas involved 
 in their work.

 Let us begin by explaining the main question in  a little more detail.  The 
 number of primes up to $x$, denoted by $\pi(x)$, is roughly $x/\log x$ for large values of 
 $x$;  this is the celebrated Prime Number Theorem\footnote{Here, and throughout, $\log$ stands for 
 the natural logarithm.}.  Therefore, if we randomly choose an integer 
 near $x$,  then it has about a $1$ in $\log x$ 
 chance of being prime.  
 In  other words, as we look at primes around size $x$, the average gap between 
 consecutive primes is about $\log x$.  As $x$ increases, the primes 
 get sparser, and the gap between consecutive primes tends to increase.  Here are some 
 natural questions about these gaps between prime numbers.  
 Do the gaps always remain roughly about size $\log x$, 
 or do we sometimes get unexpectedly large gaps and sometimes surprisingly 
 small gaps? Can we say something about the statistical distribution of these gaps?  That is, 
 can we quantify how often the gap is between, say, $\alpha \log x$ and $\beta \log x$, 
 given $0\le \alpha <\beta$?  Except for the primes $2$ and $3$, clearly the gap between consecutive   primes must be even. 
 Does every even number occur infinitely often as a gap between consecutive 
 primes?  For example, the twin prime conjecture says that the gap $2$ occurs infinitely. 
 How frequently should we expect the occurrence of twin primes?  
 
 Number theorists believe they know the answers to all these questions, 
 but cannot always prove that the answers are correct.  Before discussing 
 the answers let us address a possible meta-question.  Problems 
 like twin primes, and the Goldbach conjecture involve adding and subtracting 
 primes.  The reader may well wonder if such questions are natural, or 
 just isolated curiosities.  After all, shouldn't we be multiplying with primes 
 rather than adding/subtracting them? There are several possible responses to this 
 objection.  
 
 Firstly, many number theorists and 
 mathematical physicists are interested in understanding spacing statistics of various sequences of 
 numbers occurring in nature.  Examples of such sequences are prime numbers, 
 the ordinates of zeros of the Riemann zeta-function (see [21] and [23]), energy levels of large 
 nuclei, the fractional parts of $\sqrt{n}$ for $n\le N$ (see [7]), etc.  Do the spacings 
 behave like the gaps between randomly chosen numbers, or do they follow 
 more esoteric laws?  Our questions on gaps between primes fit naturally 
 into this framework.     
 
 Secondly, many additive questions on primes have applications to 
 other problems in number theory.  For example, consider primes $p$ 
 for which $2p+1$ is also a prime.  Analogously to twin primes, 
 it is conjectured that there are infinitely many such prime pairs 
 $p$ and $2p+1$.  Sophie Germain came up with these pairs in her 
 work on Fermat's last theorem.  If there are infinitely many 
 Germain pairs $p$ and $2p+1$ with $p$ lying in 
 a prescribed arithmetic progression, then Artin's primitive root conjecture --- every 
 positive number $a$ which is not a perfect square is a primitive root\footnote{That is, 
 $a$ generates the multiplicative group of residues modulo that prime.} 
 for infinitely many primes ---  would follow.  For example, if $p$ lies in the
 progression $3 \pmod {40}$, and $2p+1$ is 
 prime, then $10$ is a primitive root modulo $ 2p+1$, and 
 as Gauss noticed (and the reader can check) the decimal expansion of $1/(2p+1)$ 
 has exactly $2p$ digits that repeat.  There are also connections between 
 additive questions on primes and zeros of the Riemann zeta and other related functions.  
 Precise knowledge 
 of the frequency with which prime pairs $p$ and $p+2k$ occur (for an even number 
 $2k$) has subtle implications for the distribution of spacings between 
 ordinates of zeros of the Riemann zeta-function (see [1] and [23]).  Conversely, weird (and unlikely)
 patterns in zeros of zeta-like functions would imply the 
 existence of infinitely many twin primes (see [17])!  
 
 Finally, 
 these `additive' questions on primes are lots of fun, have led to much beautiful mathematics, and 
 inspired many generations of number theorists!
 
 \vskip .1 in 
 
 \noindent{\bf Cram{\' e}r's model.}  A useful way to think about  statistical questions 
 on prime numbers is the random --- also known as Cram{\' e}r ---  model.  The principle, based 
 on the fact that a number of size about $n$ has a $1$ in $\log n$ chance of 
 being prime, is this: 
 
 \smallskip
 {\sl The indicator function for the set of primes (that is, the 
 function whose value at $n$ is $1$ or $0$ depending on whether $n$ is prime 
 or not) behaves roughly like a sequence of independent, Bernoulli random variables $X(n)$ 
 with parameters $1/\log n$ ($n\ge 3$).   In other words, for $n \ge 3$, the random 
 variable $X(n)$ takes the value $1$ ($n$ is `prime') with probability $1/\log n$, 
 and $X(n)$ takes the value $0$ ($n$ is `composite') with probability $1-1/\log n$. For completeness, 
 let us set $X(1)=0$, and $X(2)=1$.} 
 \smallskip
 This must be taken with a liberal dose of salt: a number is either prime or composite, 
 probability does not enter the picture!  Nevertheless, the Cram{\' e}r model 
 is very effective in predicting answers, although it does have its limitations (for example, 
 if $n>2$ is prime then certainly $n+1$ is not, so the events of $n$ and $n+1$ being 
 prime are clearly not independent)  and 
 sometimes leads to incorrect predictions.  
 
 Let us use the Cram{\' e}r model to predict the probability that, given a large
 prime $p$, the next prime lies somewhere between $p +\alpha \log p$ and 
 $p+\beta \log p$.  In the Cram{\'e}r model, let $p$ be large and suppose that $X(p)=1$.  
 What is the probability that $X(p+1)=X(p+2) =\ldots =X(p+h-1)=0$ and 
 $X(p+h)=1$, for some integer $h$ in the interval $[\alpha \log p, \beta \log p]$? 
 We will find this by calculating the desired probability for a given $h$ in that interval, 
 and summing that answer over all such $h$.  For a given $h$ the probability we 
 seek is 
 $$
  %
\Big( 1-\frac{1}{\log (p+1)} \Big) \Big(1-\frac{1}{\log (p+2)}\Big) 
\cdots \Big(1-\frac{1}{\log (p+h-1)}\Big) \frac{1}{\log (p+h)}. 
$$
Since $p$ is large,  and $h$ is small compared to $p$ (it's only of size 
about $\log p$) we estimate that $\log (p+j)$ is very nearly $\log p$ for $j$ 
between $1$ and $h$.  Therefore our probability above is 
approximately $(1-1/\log p)^{h-1} (1/\log p)$, and since $1-1/\log p$ is 
about $e^{-1/\log p}$, this is roughly 
$$
e^{-(h-1)/\log p}\Big( \frac{1}{\log p}\Big).
$$  
Summing over the appropriate $h$, we find that 
the random model prediction for the probability that 
the next prime larger than $p$ lies in $[p+\alpha \log p, p+\beta\log p]$ 
is  
$$
\sum_{\alpha \log p \le h\le \beta \log p } e^{-(h-1)/\log p} \frac{1}{\log p} 
\approx \int_\alpha^{\beta} e^{-t} dt,
$$
since the left hand side looks like a Riemann sum approximation to  the integral.  

\proclaim{Conjecture 1}  Given an interval 
$0\le \alpha < \beta$, as $x\to \infty$ we have 
$$
\frac{1}{\pi(x)} \#\{ p \le x :  p_{\text{next}} \in (p +\alpha \log p, p +\beta \log p) \} 
\to \int_{\alpha}^{\beta} e^{-t} dt,
$$
where $p_{\text{next}}$ denotes the next prime larger than $p$. Here, and 
throughout the paper, the letter $p$ is reserved for primes.
\endproclaim

We have deliberately left the integral unevaluated, to suggest that 
there is a probability density $e^{-t}$ of finding $(p_{\text{next}}-p)/\log p$ 
close to $t$.  If we pick $N$ random numbers uniformly and independently 
from the interval $[0,N]$, and arrange them in ascending order, 
then, almost surely, the consecutive spacings have the probability 
density $e^{-t}$.  Thus, the Cram{\' e}r model indicates that 
the gaps between consecutive primes are distributed like the 
gaps between about $x/\log x$ numbers chosen uniformly and independently 
from the interval $[0,x]$.  In probability terminology, this 
is an example of what is known as a `Poisson process.'  

There are several related predictions we could make using the random model.  
For example, choose a random number $n$ below $x$, 
and consider the interval $[n,n+\log n]$.  The expected 
number of primes in such an interval is about $1$, by the prime number 
theorem. But of course some intervals may contain 
no prime at all while others may contain several primes.  Given a non-negative number $k$, what is 
the probability that such an interval contains exactly $k$ primes?  The reader 
may enjoy the pleasant calculation which predicts that, for large $x$, the answer 
is nearly $\frac{1^k}{k!} e^{-1}$ --- 
the answer is written so as to suggest a Poisson distribution with parameter $1$.
  
 Conjecture 1 makes clear that 
 there is substantial variation in the gaps between consecutive primes.    
 Given any large number $\Lambda$ we expect that with probability about 
 $e^{-\Lambda}$ (a tiny, but positive probability), the gap between consecutive primes is more than 
 $\Lambda$ times the average gap.  Given any small positive 
 number $\epsilon$ we expect that with probability about $1-e^{-\epsilon}$ (a small, but 
 positive probability),  the gap between consecutive primes is at most $\epsilon$ times the 
 usual gap.  Thus, two consequences of Conjecture 1 are 
 $$
 \limsup_{p\to \infty} \frac{p_{\text{next}}-p}{\log p} 
 =\infty,
 $$
 and 
 $$
 \liminf_{p\to \infty} \frac{p_{\text{next}} -p}{\log p} =0.
 $$

 \vskip .1 in 
 
 \noindent {\bf Large gaps.}  Everyone knows how to 
 construct arbitrarily long intervals of composite numbers:  just look at 
 $m!+2$, $m! +3$, $\ldots$, $m!+m$ for any natural number $m\ge 2$.   
 This shows that $\limsup_{p\to \infty} (p_{\text{next}}-p) = \infty$.  However, if 
 we think of $m!$ being of size about $x$ then a little calculation 
 with Stirling's formula shows that $m$ is about size $(\log x )/\log \log x$.
 We realize, with dismay, that the `long' gap we have constructed is not even as 
 large as the average gap of $\log x$ given by the prime number theorem.  
 A better strategy is to take $N$ to be the product of 
 the primes that are at most $m$, and note again that $N+2$, 
 $\ldots$, $N+m$ must all be composite.   It can be shown that $N$ is 
 roughly of size $e^m$.  Thus we have found a 
 gap at least about $\log N$, which is better than before, but still 
 not better than average.  Can we modify the argument a little?  In creating 
 our string of $m-1$ consecutive composite numbers, we forced these 
 numbers to be divisible by some prime below $m$.  Can we somehow use 
 primes larger than $m$ to force $N+m+1$, $N+m+2$, etc., to be composite, and 
 thus create longer chains of composite numbers?  
 In the 1930s, in a series of papers Westzynthius [27], Erd{\H o}s [8] and 
 Rankin [25] found ingenious ways of making this idea work.  
 The best estimate was obtained by Rankin, who 
 proved that there exists a positive constant $c$ 
 such that for infinitely many primes $p$, 
 $$
 p_{\text{next}} -p > c \log p \frac{(\log \log p) \log \log \log \log p}{(\log \log \log p)^2}.
  $$ 
 The fraction above does grow\footnote{although so slowly that, as the joke goes, no one 
 has observed it doing so!}, and so 
 $$
 \limsup_{p\to \infty} \frac{p_{\text{next}}-p}{\log p} = \infty,
 $$
 as desired.  We 
 should remark here that, although very interesting work has been done on improving the constant $c$ above,  
 Rankin's result provides the largest known gap between primes.  Erd{\H o}s 
 offered \$10,000 for a similar conclusion involving a faster growing function.  
 Bounty hunters may note that the largest Erd{\H o}s prize that has 
 been collected is \$1,000, by Szemeredi [26] for his marvellous result on 
 the existence of long arithmetic progressions in sets of positive density.
 
 What should we conjecture 
 for the longest gap between primes? Cram{\' e}r's 
 model suggests that 
 $$
 \limsup_{p \to \infty} \frac{p_{\text{next}}-p}{(\log p)^2} = c, \tag{1}
 $$
 with $c=1$.  
The rationale behind this is that the probability that $X(n)=1$ and that the next `prime' is bigger than  
 $n+ (1+\epsilon) \log^2 n$ is about $1/(n^{1+\epsilon}\log n)$, 
 by a calculation similar to the one leading up to Conjecture 1.   
 If $\epsilon$ is negative the sum of this 
 probability over all $n$ diverges, and the Borel-Cantelli lemma tells us that, almost surely, 
 such long gaps  occur infinitely often.  If $\epsilon$  is positive, the corresponding sum converges 
 and the Borel-Cantelli lemma says that almost surely we get these longer gaps only a 
 finite number of times.  More sophisticated analysis has however revealed that 
 (1) is one of those questions which expose the limitations of the Cram{\' e}r model.  
 It appears unlikely that the value of $c$ is $1$ as predicted 
 by the Cram{\' e}r model, and that $c$ should be at least $2e^{-\gamma}\approx 1.1229$ where 
$\gamma$ is Euler's constant.  No one has felt brave enough to suggest what the precise 
value of $c$ should be!  This is because (1) is far beyond what `reasonable' conjectures 
such as the Riemann hypothesis would imply.  An old conjecture says 
that there is always a prime between two consecutive squares.  Even this lies 
(slightly) beyond the reach of the Riemann hypothesis, and all it would imply is that 
$$
\limsup_{p\to \infty} \frac{p_{\text{next}}-p}{\sqrt{p}} \le 4; 
$$
a statement much weaker than (1) with a finite value of $c$. 

We cut short our discussion on long gaps here, since our focus will be on small gaps; for more information on these and 
related problems, we refer the reader to the excellent survey articles by Heath-Brown [18] and Granville [15].

 \vskip .1 in 
 \noindent {\bf Small gaps.}  Since the average spacing between $p$ and 
 $p_{\text{next}}$ is about $\log p$, clearly 
 $$
 \liminf_{p\to \infty} \frac{p_{\text{next}}-p}{\log p} \le 1.
 $$
 Erd{\H o}s [9] was the first to show that the $\liminf$ is strictly less than $1$.  
 Other landmark results in the area are 
 the works of Bombieri and Davenport [3], Huxley [20], and Maier [22], who introduced 
 several new ideas to this study and progressively reduced the $\liminf$ to $\le 0.24\ldots$.  
 Enter Goldston, Pintz, and \Yildirim\!\!:
 
 \proclaim{Theorem 1} We have 
 $$
 \liminf_{p\to \infty} \frac{p_{\text{next}}-p}{\log p} = 0.
 $$
 \endproclaim

So there are substantially smaller gaps between primes than the average!  
What about even smaller gaps?  Can we show that $\liminf_{p\to \infty} 
(p_{\text{next}} -p ) < \infty$ (bounded gaps), or perhaps even 
$\liminf_{p\to \infty} (p_{\text{next}}-p) =2$ (twin primes!)?   

\proclaim{Theorem 2}  Suppose the Elliott-Halberstam conjecture on the 
distribution of primes in arithmetic progressions holds true. 
  Then 
$$
\liminf_{p \to \infty} (p_{\text{next}}-p) \le 16.
$$
\endproclaim 

What is the Elliott-Halberstam conjecture? One valuable 
thing that we know about primes is their distribution in 
arithmetic progressions.  Knowledge of this, in the 
form of the Bombieri-Vinogradov theorem, plays a crucial 
role in the proof of Theorem 1.  To obtain the stronger 
conclusion of Theorem 2, one needs a better understanding 
of the distribution of primes in progressions and the 
Elliott-Halberstam conjecture provides the necessary stronger 
input.  Vaguely, the Goldston-Pintz-\Yildirim results say 
that if the primes are well separated with no small gaps between them,  
then something weird must happen to their distribution 
in progressions.  

Given a progression $a\pmod q$ let $\pi(x;q,a)$ denote the 
number of primes below $x$ lying in this progression.  
Naturally we may suppose that $a$ and $q$ are coprime, else 
there is at most one prime in the progression. Now there are 
$\phi(q)$ --- this is Euler's $\phi$-function --- such 
progressions $a \pmod q$ with $a$ coprime to $q$.  We would expect 
that each progression captures its fair share of primes.  In other 
words we expect that $\pi(x;q,a)$ is roughly $\pi(x)/\phi(q)$.  
The prime number theorem in arithmetic progressions tells us that this is 
true if we view $q$ as being fixed and let $x$ go to infinity.  

In applications, such as Theorem 1, we need information on 
$\pi(x;q,a)$ when $q$ is not fixed, but growing with $x$.  When 
$q$ is growing slowly, say $q$ is like $\log x$, the prime number 
theorem in arithmetic progressions still applies.  However if 
$q$ is a little larger, say $q$ is of size $x^{\frac 13}$, 
then currently we cannot prove the equidistribution of primes 
in the available residue classes $\pmod q$.  Such a result would 
be implied by the Generalized Riemann Hypothesis (indeed for $q$ 
up to about $\sqrt{x}$), but of course the Generalized Riemann 
Hypothesis remains unresolved.  In this context, Bombieri and Vinogradov 
showed that the equidistribution of primes in progressions holds, 
not for each individual $q$, but on average over $q$ (that is, for 
a typical $q$) for $q$ going up to about $\sqrt{x}$.  Their result 
may be thought of as the `Generalized Riemann Hypothesis on average.'
  
The Elliott-Halberstam conjecture says that the equidistribution of 
primes in progressions continues to hold on average for $q$ going 
up to $x^{1-\epsilon}$ for any given positive $\epsilon$.  In some 
ways, this lies deeper than the Generalized Riemann Hypothesis which
permits only $q\le \sqrt{x}$.   

We hope that the reader has formed a rough impression of 
the nature of the assumption in Theorem 2.  We will state the 
Bombieri-Vinogradov theorem and Elliott-Halberstam conjecture 
precisely in the penultimate section devoted to primes 
in progressions.  

\vskip .1 in 

\noindent{\bf The Hardy-Littlewood conjectures.} We already noticed a 
faulty feature of the 
Cram{\' e}r model:  given a large prime $p$, the probability that $p+1$ is prime is not $1/\log (p+1)$ 
but $0$ because $p+1$ is even.  
Neither would we expect the conditional probability of $p+2$ being prime to 
be simply $1/\log (p+2)$: after all, $p+2$ is guaranteed to be odd and this should give it 
a better chance of being prime.   How should we formulate the correct probability for 
$p+2$ being prime?  More precisely, what should be the conjectural asymptotics 
for 
$$
\# \{ p\le x: \ p+2 \text{ prime} \} ?
$$
The Cram{\' e}r model would have predicted that this is about $x/(\log x)^2$. 
While we must definitely modify this, it also seems reasonable that $x/(\log x)^2$ is 
the right size for the answer.  So maybe the answer is about $cx/(\log x)^2$ for an 
appropriate constant $c$.  

Long ago Hardy and Littlewood [16] figured out what the 
right conjecture should be.  The problem with the Cram{\' e}r model is 
that it treats $n$ and $n+2$ as being independent, whereas they are clearly dependent.  
If we want $n$ and $n+2$ both to be prime, then they must both 
be odd, neither of them must be divisible by $3$, nor by $5$, and so on.  
If we choose $n$ randomly, the probability that $n$ and $n+2$ are both odd is $1/2$.  
In contrast, two randomly chosen numbers would both be odd with a $1/4$ probability.
If neither $n$ nor $n+2$ is divisible by $3$ then $n$ must be $2 \pmod 3$, which 
has a $1/3$ probability.   
On the other hand, the probability that two randomly chosen numbers 
are not divisible by $3$ is $(2/3)\cdot (2/3)= 4/9$. 
Similarly, for any prime $\ell \ge 3$, 
 the probability that $n$ and $n+2$ are not divisible by $\ell$ 
is $1-2/\ell$, which is a little different from the probability $(1-1/\ell)^2$ that two randomly chosen 
integers are both not divisible by $\ell$.  For the prime $2$ we 
must correct the probability $1/4$ by multiplying by  $2=(1-1/2)(1-1/2)^{-2}$, 
and for all primes $\ell \ge 3$ we must correct the probability $(1-1/\ell)^2$ 
by multiplying by $(1-2/\ell)(1-1/\ell)^{-2}$.  The idea is that if we multiply all 
these correction factors together then we have accounted for `all the ways' 
in which $n$ and $n+2$ are dependent, producing the required correction constant $c$.  
Thus the conjectured value for $c$ is the product over primes 
$$
\Big(1-\frac 12\Big) \Big( 1-\frac 1{2}\Big)^{-2} \prod_{\ell \ge 3} 
\Big(1-\frac 2\ell\Big) \Big(1-\frac 1\ell\Big)^{-2}.
$$

Let us make a synthesis of the argument above, which will allow us 
to generalize it.  For any prime $\ell$ let $\nu_{\{0,2\}}(\ell)$ 
denote the number of distinct residue classes $(\bmod \ \ell)$ occupied by the numbers $0$ 
and $2$.  If we want $n$ and $n+2$ to be both coprime to $\ell$ then $n$ 
must $n$ must avoid the residue classes occupied by $-0$ and $-2 \pmod \ell$, 
so that $n$ must lie in one of $\ell - \nu_{\{0,2\}}(\ell)$ residue classes.  The probability 
that this happens is $1-\nu_{\{0,2\}}(\ell)/\ell$, so the correction 
factor for $\ell$ is $(1-\nu_{\{0,2\}}(\ell)/\ell)(1-1/\ell)^{-2}$.  
As before, consider the infinite product over primes
 $$
 {\frak S}(\{0,2\}) := \prod_{\ell} \Big( 1-\frac{\nu_{\{0,2\}} (\ell)}{\ell}\Big) \Big(1-\frac{1}{\ell}\Big)^{-2}.   
 $$
 The infinite product certainly converges: the terms for $\ell \ge 3$ are all less 
 than $1$ in size.  Moreover, it converges to a non-zero number.  Note that none of the 
 factors above is zero, and that for large $\ell$ the logarithm of the corresponding 
 factor above is very small --- it is $\log (1-1/(\ell -1)^2) \approx - 1/\ell^2$.  Thus  
 the sum of the logarithms converges, and the product is non-zero; indeed ${\frak S}(\{0,2\})$ 
 is numerically about $1.3203$.   
 Then the conjecture is that for large $x$ 
 $$
 \# \{ p\le x: \ p+2 \text{ prime} \} \sim {\frak S}(\{0,2\}) \frac{x}{(\log x)^2}.
 $$
 Here and below, the notation $f(x) \sim g(x)$ means that $\lim_{x\to \infty} f(x)/g(x)=1$.

 The conjecture generalizes readily:  Suppose we are given a set ${\Cal H}=\{h_1, h_2, \ldots, h_k\}$ 
 of non-negative integers and we want to find the frequency with which $n+h_1$, $\ldots$, $n+h_k$ are all 
 prime.  For a prime number $\ell$, we define $\nu_{{\Cal H}}(\ell)$ to be the 
 number of distinct residue classes $(\bmod \ \ell)$ occupied by ${\Cal H}$.  We define 
 the `singular series'\footnote{The terminology is not entirely whimsical:  Hardy and Littlewood 
 originally arrived at their conjecture through a heuristic application of their `circle method.'  
 In their derivation, ${\frak S}({\Cal H})$ did arise as a series rather than as our product.} 
 $$
 {\frak S}({\Cal H}) = \prod_{\ell} \Big(1-\frac{\nu_{\Cal H}(\ell)}{\ell}\Big) \Big( 1-\frac 1\ell\Big)^{-k}. 
 \tag{2}
 $$
 If $\ell$ is larger than all elements of ${\Cal H}$ then $\nu_{\Cal H}(\ell)=k$, 
 and for such $\ell$ the terms in the product are less than $1$.  Thus 
 the product converges.  When does it converge to a non-zero number?   If $\nu_{\Cal H}(\ell) =
 \ell$ for some prime $\ell$ then one of the terms in our product vanishes, and so 
 our product must be zero.  Suppose none of the terms is zero.  For 
 large $\ell$ the logarithm of the corresponding factor is 
$$
\log \Big(1-\frac k\ell\Big)\Big(1-\frac 1\ell\Big)^{-k} \approx - \frac{k(k+1)}{2\ell^2},
$$
and so the sum of the logarithms converges, and our product is non-zero.  Thus the 
singular series is zero if and only if $\nu_{\Cal H}(\ell) =\ell$ for some prime $\ell$ 
--- that is, if and only if the numbers $h_1, \ldots, h_k$ occupy {\sl all} 
the residue classes $(\bmod \ \ell)$ for some prime $\ell$.  In that case, for any 
$n$ one of the numbers $n+h_1$, $\ldots$, $n+h_k$ must be a multiple of $\ell$, 
and so there are only finitely many prime $k$-tuples $n+h_1$, $\ldots$, $n+h_k$.

\proclaim{The Hardy-Littlewood conjecture}   Let ${\Cal H} 
= \{h_1, \ldots, h_k\}$ be a set of positive integers such that 
${\frak S}({\Cal H}) \neq 0$.  Then 
 $$
 \# \{ n\le x: \ n+h_1, \ldots, n+h_k \text{ prime} \} 
 \sim {\frak S}({\Cal H}) \frac{x}{(\log x)^k}. 
 $$
 \endproclaim
 
 
 It is easy to see that ${\frak S}(\{0,2r\}) \neq 0$ for every non-zero even number $2r$.  Thus the 
 Hardy-Littlewood conjecture predicts that there are about ${\frak S}(\{0,2r\}) x/(\log x)^2$ 
 prime pairs $p$ and $p+2r$ with $p$ below $x$.  Further, the number of these pairs for which 
 $p+2d$ is prime for some $d$ between $1$ and $r-1$ is at most a constant times $x/(\log x)^3$.  
 We deduce that there should be infinitely many primes $p$ for which 
 the gap to the next prime is exactly $2r$.   Thus every positive even number should occur infinitely often 
 as a gap between successive primes, but we don't know this for a single even number! 
 
 For any $k$, it is easy to find $k$-element sets ${\Cal H}$ with ${\frak S}({\Cal H}) \neq 0$. 
 For example, take ${\Cal H}$ to be any $k$ primes all larger than $k$.  Clearly if 
 $\ell>k$ then $\nu_{{\Cal H}}(\ell) \le k <\ell$, while if $\ell \le k$ then the residue class 
 $0 \pmod \ell$ must be omitted by the elements of ${\Cal H}$ (they are primes!) and 
 so once again $\nu_{\Cal H}(\ell) <\ell$.

 We make one final comment before turning (at last!) to the ideas behind the proofs of 
 Theorems 1 and 2.  Conjecture 1 was made on the strength of the Cram{\' e}r model, but we 
 have just been discussing how to modify the Cram{\'e}r probabilities for prime $k$-tuples. 
 A natural question is whether the Hardy-Littlewood conjectures are consistent with Conjecture 1.   In a 
 beautiful calculation [11], Gallagher showed that Conjecture 1 can in fact be obtained 
 starting from the Hardy-Littlewood conjectures.  The crucial point in his proof 
 is that although ${\frak S}({\Cal H})$ is not always $1$ (as the Cram{\' e}r model would have), 
 it is {\sl approximately} $1$ on average over all $k$-element sets ${\Cal H}$ with the $h_j \le h$.  That is, 
 as $h \to \infty$,  
 $$
 \sum\Sb1\le  h_1 < h_2 < \ldots < h_k \le h \endSb {\frak S}(\{h_1,\ldots,h_k\}) 
 \sim \sum\Sb 1\le h_1 < h_2 <\ldots < h_k \le h \endSb 1. \tag{3}
 $$
 
 \vskip .1 in

 \noindent {\bf The ideas of Goldston, Pintz and \Yildirim\!\!. }  We will 
 start with the idea behind Theorem 2.  Let $k$ be a given positive integer 
 which is at least $2$.  Let ${\Cal H} = \{h_1 < \ldots< h_k\}$ be a set 
 with ${\frak S}({\Cal H}) \neq 0$.  We aspire to the Hardy-Littlewood conjecture 
 which says that there must be infinitely many $n$ such that $n+h_1$, $\ldots$, $n+h_k$ 
 are all prime.  Since there are infinitely many primes, trivially at least one of the 
 numbers $n+h_1$, $\ldots$, $n+h_k$ is prime infinitely often.  
 Can we do a little better: can we show that two of the numbers $n+h_1$, 
 $\ldots$, $n+h_k$ are prime infinitely often?  If we could, then we 
 would plainly have that $\liminf_{p\to \infty} (p_{\text{next}} -p) \le (h_k-h_1)$.

 How do we detect two primes in $n+h_1$, $\ldots$, $n+h_k$?  
 Let $x$ be large and consider $n$ varying between $x$ and $2x$.  
 Suppose we are able to find a function $a(n)$ which is always non-negative, and 
 such that, for each $j=1$, $\ldots$, $k$,  
 $$
 \sum\Sb x\le n \le 2x \\ n+h_j \text{ prime}\endSb a(n) 
 > \frac{1}{k} \sum\Sb x\le n\le 2x \endSb a(n). \tag{4} 
 $$
 Then summing over $j=1$, $\ldots$, $k$,  it would follow that 
 $$
 \sum\Sb x\le n\le 2x \endSb \#\{ 1\le j\le k: n+h_j \text{ prime}\}\  a(n) 
> \sum\Sb x\le n\le 2x \endSb a(n),
$$
so that for some number $n$ lying between $x$ and $2x$ 
we must have at least two primes among 
$n+h_1$, $\ldots$, $n+h_k$.  

 Of course, the question is how do we find such a function $a(n)$ satisfying (4)?  
 We would like to take $a(n)=1$ if $n+h_1$, $\ldots$, $n+h_k$ are 
 all prime, and $0$ otherwise.  But then the problem 
 of evaluating $\sum_{x\le n\le 2x} a(n)$ is precisely that 
 of establishing the  Hardy-Littlewood conjecture.

 The answer is suggested by sieve theory, especially the theory of Selberg's sieve.  
 Sieve theory is concerned with finding primes, or numbers without too many 
 prime factors, among various integer sequences.  Some of 
 the spectacular achievements of this theory are Chen's theorem [5] that 
 for infinitely many primes $p$, the number $p+2$ has at most two prime factors; the 
 result of Friedlander and Iwaniec [10] that there are infinitely many primes of 
 the form $x^2+y^4$ where $x$ and $y$ are integers; 
 and the result of Heath-Brown [19] that there are infinitely many 
 primes of the form $x^3 +2y^3$ where $x$ and $y$ are integers.  We recall here very briefly the idea 
 behind Selberg's sieve.  
 
 \vskip .1 in 
 
 \noindent{\bf Interlude on Selberg's sieve.}  We illustrate Selberg's sieve by giving an 
 upper bound on the number of prime $k$-tuples $n+h_1$, $\ldots$,  $n+h_k$ with $x\le n\le 2x$.  
 The idea is to find a `nice' function $a(n)$ which equals $1$ 
 if $n+h_1$, $\ldots$, $n+h_k$ are all prime, and is non-negative otherwise.  
 Then $\sum_{x\le n\le 2x} a(n)$ 
 provides an upper bound for the number of prime $k$-tuples.  Of course, we 
 must choose $a(n)$ appropriately, so as to be able to evaluate $\sum_{x\le n\le 2x}a(n)$.

 Selberg's choice for $a(n)$ is as follows:  Let $\lam_d$ be a sequence of real numbers 
 such that
 $$
 \lam_1 =1, \qquad \text{and with } \qquad \lam_d =0 \qquad \text{for } d>R. \tag{5}
 $$
 Choose\footnote{Below, the symbol $a|b$ means that $a$ divides $b$.}
 $$
 a(n) =\Big(\sum_{d|(n+h_1)\ldots (n+h_k)}  \lam_d\Big)^2. \tag{6}
 $$
 Being a square, $a(n)$ is clearly non-negative.  If $R <x\le n$ 
 and $n+h_1$, $\ldots$, $n+h_k$ are all prime, then the only non-zero term in (6) is for $d=1$ and so $a(n)=1$ as desired.   Therefore we assume that $R< x $ below.  The goal is to choose $\lam_d$ so as to minimize 
 $\sum_{x\le n\le 2x}a(n)$.  There is an advantage to allowing $R$ as large as possible, 
 since this gives us greater flexibility in choosing the parameters $\lam_d$.  On 
 the other hand it is easier to estimate $\sum_{x\le n\le 2x} a(n)$ when $R$ is 
 small since there are fewer divisors $d$ to consider.  In the problem at hand, it turns out 
 that we can choose $R$ roughly of size  $\sqrt{x}$.  This choice leads to an upper bound 
 for the number of prime $k$-tuples of about $2^k \cdot k! {\frak S}({\Cal H}) x/(\log x)^k$.  
 That is, a bound about $2^k \cdot k!$ times the conjectured Hardy-Littlewood 
 asymptotic.  
 
 Expanding out the square in (6) and summing over $n$, we must evaluate 
 $$
 \sum\Sb d_1, d_2 \endSb \lam_{d_1}\lam_{d_2} 
 \sum\Sb x\le n\le 2x\\ d_1 | (n+h_1)\cdots (n+h_k) \\ d_2 |(n+h_1)\cdots (n+h_k) \endSb 1 
 = \sum\Sb d_1, d_2 \endSb \lam_{d_1}\lam_{d_2} 
 \sum\Sb x\le n\le 2x\\ [d_1,d_2] | (n+h_1)\cdots (n+h_k) \endSb 1,
 $$
 where $[d_1,d_2]$ denotes the l.c.m. of $d_1$ and $d_2$.  The condition $[d_1,d_2]|(n+h_1)\cdots 
 (n+h_k)$ means that $n$ must lie in a certain number (say, $f([d_1,d_2])$) 
 of residue classes $\pmod{[d_1,d_2]}$.  Can we count the number of $x\le n\le 2x$ lying 
 in the union of these arithmetic progressions?  Divide the interval $[x,2x]$ into intervals of 
 length $[d_1,d_2]$ with possibly one smaller interval left over at the end.  
 Each complete interval (and there are about $x/[d_1,d_2]$ of these) 
 gives $f([d_1,d_2])$ values of $n$; the last shorter interval contributes an indeterminate 
 `error' between $0$ and $f([d_1,d_2])$.  So, at least if $[d_1,d_2]$ is a bit smaller 
 than $x$, we can estimate the sum over $n$ accurately.  Since $[d_1,d_2] \le d_1d_2 \le R^2$, 
 if $R$ is a bit smaller than\footnote{To be precise, $R$ must be $\le \sqrt{x}/(\log x)^{2k}$, say.}  $\sqrt{x}$, 
 then the sum over $n$ can be evaluated accurately.  Let us suppose 
 that $R$ is about size $\sqrt{x}$ and that the error terms can be disposed of satisfactorily.   
 It remains to handle the main term contribution to $\sum_{x\le n\le 2x} a(n)$, namely
 $$
x \sum\Sb d_1,d_2\le R\endSb \frac{f([d_1,d_2])}{[d_1,d_2]} \lam_{d_1}\lam_{d_2}. \tag{7}
$$
The reader may wonder what $f([d_1,d_2])$ is. Let us work this out in the case 
when $[d_1,d_2]$ is not divisible by the square of any prime; the other case is more 
complicated, but not very important in this problem.  If $p$ is a prime and 
we want $p|(n+h_1)\cdots(n+h_k)$ then clearly $n \equiv -h_j \pmod p$ 
for some $j$,  so that $n$ lies in one of $\nu_{\Cal H}(p)$ residue classes $\pmod p$.  
By the chinese remainder theorem it follows that if $[d_1,d_2]|(n+h_1)\cdots (n+h_k)$ 
then $n$ lies in $\prod_{p|[d_1,d_2]} \nu_{\Cal H}(p)$ residue classes $\pmod {[d_1,d_2]}$.  
Thus $f$ is a multiplicative function\footnote{These are functions satisfying $f(mn)=f(m)f(n)$ 
for any pair of coprime integers $m$ and $n$.}, with $f(p) =\nu_{\Cal H}(p)$.

The problem in Selberg's sieve is to choose $\lam_d$ subject to the linear constraint (5) 
in such a way as to minimize the quadratic form (7) (that would give the best upper 
bound for $\sum_{x\le n\le 2x} a(n)$).  This can be achieved using Lagrange multipliers, 
or by diagonalizing the quadratic form (7).   We do not give the details of this 
calculation but just record the result  obtained.  The optimal choice of $\lam_d$ 
for $d\le R$ is given by 
$$
\lam_d \approx \mu(d) \Big(\frac{\log R/d}{\log R}\Big)^k,
$$
where $\mu(d)$ is the M{\" o}bius function.\footnote{$\mu(d)=0$ if $d$ is divisible 
by the square of a prime.  Otherwise $\mu(d)= (-1)^{\omega(d)}$ 
where $\omega(d)$ is the number of distinct primes dividing $d$.}   With this choice of $\lam_d$ 
the quantity in (7) is 
$$
\approx k! {\frak S}({\Cal H}) \frac{x}{(\log R)^k} \approx 2^k \cdot k! {\frak S}({\Cal H}) 
\frac{x}{(\log x)^k}.
$$
 
 The appearance at this stage of the M{\" o}bius function is not surprising, as 
 it is very intimately connected with primes.  For example, the reader can check that 
 $\sum_{d|m} \mu(d) (\log m/d )^k$ equals $0$ unless $m$ is divisible by 
 at most $k$ distinct prime factors.  When $m=p_1\cdots p_k$ is the 
 product of $k$ distinct prime factors it equals $k! (\log p_1)
 \cdots (\log p_k)$, and there is a more complicated formula if $m$ is composed 
 of fewer than $k$ primes, or if $m$ is divisible by powers of primes. 
  Applying this to $m=(n+h_1) \cdots (n+h_k)$, we are 
 essentially picking out prime $k$-tuples!  The optimum in Selberg's sieve is a kind 
 of approximation to this identity.  
 
 \vskip .1 in

 \noindent{\bf Return to Goldston-Pintz-\Yildirim\!\!.}  We want to find a non-negative function 
 $a(n)$ so as to make (4) hold.  Motivated by Selberg's sieve we may try to 
find optimal $\lam_d$ as in (5) and again choose $a(n)$ as in (6).  If we try such a 
choice, then our problem now is to maximize the ratio 
$$
 \Big(\sum\Sb x\le n\le 2x\\ n+h_j \text{prime} \endSb a(n) \Big) \Big/  \Big(
\sum_{x\le n\le 2x} a(n)\Big). \tag{8}
$$
We'd like this ratio to be $>1/k$.  Notice again that it is advantageous to 
choose $R$ as large as possible to give greatest freedom in choosing $\lam_d$, 
but in order to evaluate the sums above there may be restrictions on the size of $R$.  
In dealing with the denominator we saw that there is a restriction $R \le \sqrt{x}$ 
(essentially) and that in this situation the denominator in (8) is given by 
the quadratic form (7).  We will see below that in dealing with the numerator of (8), 
a more stringent restriction on $R$ must be made:  we can only take $R$ around 
size $x^{\frac 14}$.  

In any case, (8) is the ratio of two quadratic forms, and this ratio needs 
to be maximized keeping in mind the linear constraint (5).  This optimization problem 
is more delicate than the one in Selberg's sieve.  It is not clear how to proceed 
most generally: Lagrange multipliers become quite messy, and we can't quite 
diagonalize both quadratic forms simultaneously.  It helps to narrow the search to 
a special class of $\lam_d$.  Motivated by Selberg's sieve we will search for the 
optimum among the choices (for $d\le R$) 
$$
\lam_d = \mu(d) P\Big(\frac{\log R/d}{\log R}\Big).
$$
Here $P(y)$ denotes a polynomial such that $P(1)=1$ and 
such that $P$ vanishes to order at least $k$ at $y=0$.  The 
condition that $P$ be a polynomial can be relaxed a bit but this is not 
important.  It is however vital for the analysis that $P$ should vanish to 
order $k$ at $0$. Our aim is to find a choice for $P$ which makes the ratio in (8) 
large.

With this choice of $\lam_d$ we can use standard arguments 
to evaluate (7) and thus the denominator in (8).  Omitting the 
long, technical details, the answer is that for $R$ a little 
below $\sqrt x$, the denominator in (8) is 
$$
\sim \frac{x}{(\log R)^k} {\frak S}({\Cal H}) 
\int_0^1 \frac{y^{k-1}}{(k-1)!} P^{(k)}(1-y)^2 dy, \tag{9}
$$
where $P^{(k)}$ denotes the $k$-th derivative of the polynomial $P$.

To handle the numerator of (8), we expand out the square in (6) and 
sum over $x\le n\le 2x$ with $n+h_j$ being prime.  Thus the numerator is 
$$
\sum\Sb d_1, d_2 \le R \endSb 
\lam_{d_1} \lam_{d_2} \sum\Sb x\le n\le 2x \\ [d_1,d_2] | (n+h_1)\cdots (n+h_k) \\ 
n+h_j \text{prime}\endSb 1.
$$
How can we evaluate the inner sum over $n$?  As we saw before, 
the condition $[d_1,d_2]$ divides $(n+h_1)\cdots (n+h_k)$ means that $n$ lies in $f([d_1,d_2])$ 
arithmetic progressions $\pmod {[d_1,d_2]}$.  For each of these progressions 
we must count the number of $n$ such that $n+h_j$ is prime.  Of course, for 
some of the $f([d_1,d_2])$ progressions it may happen that $n+h_j$ 
automatically has a common factor with $[d_1,d_2]$ and so cannot be 
prime.  Suppose there are $g([d_1,d_2])$ progressions such that $n+h_j$ 
is guaranteed to be coprime to  $[d_1,d_2]$.  For each of these progressions 
we are counting the number of primes between $x$ and $2x$ lying 
in a reduced residue class\footnote{A reduced 
residue class $\pmod q$ is a progression $a \pmod q$ where $a$ is coprime 
to $q$.}  $\pmod {[d_1,d_2]}$.   Given a 
modulus $q$, the prime number theorem in 
arithmetic progressions says that the primes are roughly equally divided among 
the reduced residue classes $\pmod{q}$.  Thus, ignoring error terms completely, 
we expect the sum over $n$ to be about 
$$
\frac{\pi(2x) -\pi(x)}{\phi([d_1,d_2])} g([d_1,d_2]). 
$$
The $\phi([d_1,d_2])$ in the denominator is Euler's $\phi$-function: 
for any integer $m$, $\phi(m)$ counts the number of reduced residue classes $\pmod {m}$.   
Since $\pi(2x)-\pi(x)$ is about $x/\log x$ we `conclude' that the numerator in (8) is 
about 
$$
\frac{x}{\log x} \sum\Sb d_1, d_2 \le R \endSb 
\lam_{d_1} \lam_{d_2} \frac{g([d_1,d_2])}{\phi([d_1,d_2])}. \tag{10}
$$
This is the expression analogous to (7) for the numerator.

Two big questions: what is the function $g$, and for what range of $R$ can 
we handle the error terms above?  Let us first describe $g$.  As with $f$ let 
us suppose that $[d_1,d_2]$ is not divisible by the square of any prime.  As noted 
earlier, if $p$ is prime and $p|(n+h_1)\cdots (n+h_k)$ then $n$ lies in 
one of $\nu_{\Cal H}(p)$ residue classes $\pmod p$.  If we want $n+h_k$ to be 
prime, then one of these residue classes, namely $n\equiv -h_j \pmod p$, 
must be forbidden.  Thus there are now $\nu_{\Cal H}(p)-1$ residue 
classes available for $n \pmod p$.  In other words, $g(p)= \nu_{{\Cal H}}(p)-1$, 
and the chinese remainder theorem shows that $g$ must be defined multiplicatively: 
$$
g([d_1,d_2]) = \prod_{p|[d_1,d_2]} (\nu_{\Cal H}(p)-1).  
$$

We will postpone the detailed discussion on primes in arithmetic progressions which 
is needed to handle the error terms above. 
For the moment, let us note that the Bombieri-Vinogradov theorem 
(which is a powerful substitute for the generalized Riemann hypothesis
 in many applications) allows us to control $\pi(x;q,a)$ (the 
 number of primes up to $x$ which are congruent to $a\pmod q$), 
 on average over $q$, for $q$ up to about $\sqrt{x}$.   Since our moduli 
 are $[d_1,d_2]$, which go up to $R^2$, we see that $R$ may be chosen 
 up to about $x^{\frac 14}$.  Conjectures of Montgomery, and Elliott 
 and Halberstam (discussed below) would permit larger values of $R$, 
 going up to $x^{\frac 12-\epsilon}$ for any $\epsilon >0$.

 Thus, with $R$ a little below $x^{\frac 14}$, the expression (10) 
 does give a good approximation to the numerator of (8).  Now a  
 standard but technical argument can be used to evaluate (10).  
 As with (9), the answer is 
 $$
 \sim \frac{x}{(\log x) (\log R)^{k-1}} {\frak S}({\Cal H}) \int_{0}^{1} 
 \frac{y^{k-2}}{(k-2)!} P^{(k-1)}(1-y)^2 dy.
 \tag{11}
 $$

 Assuming that ${\frak S}({\Cal H})
\neq 0$, it follows from (9) and (11) that the ratio in (8) is about 
 $$
 \frac{\log R}{\log x} \Big(  \int_{0}^{1} 
 \frac{y^{k-2}}{(k-2)!} P^{(k-1)}(1-y)^2 dy\Big) \Big/ \Big( \int_{0}^{1} 
 \frac{y^{k-1}}{(k-1)!} P^{(k)}(1-y)^2 dy\Big). \tag{12} 
 $$
 This is the moment of truth: can we choose $P$ so as to make this a little 
 larger than $1/k$?  
 
 Here is a good choice for $P$: take $P(y)=y^{k+r}$ for a non-negative integer 
 $r$ to be chosen optimally.  After some calculations with beta-integrals, 
 we see that (12) then equals 
 $$
 \Big(\frac{\log R}{\log x} \Big) \Big( \frac{2(2r+1)}{(r+1) (k+2r+1)}\Big).
 $$
 This is largest when $r$ is about $\sqrt{k}/2$, and the second fraction above 
 is close to but less than $4/k$.  Since we can choose $R$ a little 
 below $x^{\frac 14}$, the first fraction is close to but less than $1/4$.  Thus 
 (12) is very close to, but less than, $1/k$.  We therefore barely fail to prove 
 bounded gaps between primes!   Of course, we just tried one choice of $P$; 
 maybe there is a better choice which gets us over the edge.  Unfortunately, 
the second fraction in (12) {\sl cannot} be made larger than $4/k$.  If we set 
$Q(y)=P^{(k-1)}(y)$ then $Q$ is a polynomial, not identically zero,  
with $Q(0)=0$; for such polynomials $Q$ we claim that the unfortunate inequality 
$$
\int_0^1 \frac{y^{k-2}}{(k-2)!} Q(1-y)^2 dy < \frac{4}{k} \int_0^1 
\frac{y^{k-1}}{(k-1)!} Q^{\prime}(1-y)^2 dy
$$
holds. The reader can try her hand at proving this.

We now have enough to prove Theorem 2!  If we can choose $R$ a little larger 
than $x^{\frac 14}$ then for suitably large $k$ the quantity in (12) can be made larger than $1/k$ as 
desired.  If we allow $R=x^{\frac 12-\epsilon}$ as the Elliott-Halberstam conjecture 
predicts, then with $k=7$ and $r=1$ we can make (12) nearly $1.05/k > 1/k$.  
Thus, if we take any set ${\Cal H}$ with seven elements and ${\frak S}({\Cal H})\neq 0$ 
then for infinitely many $n$ at least two of the numbers $n+h_1$, $\ldots$, $n+h_k$ 
are prime!  By choosing a more careful polynomial $P$ we can make do with 
six element sets ${\Cal H}$ rather than seven.  The first six primes larger than 
$6$ are $7, 11, 13, 17, 19$, and $23$, and so ${\frak S}(\{7, 11, 13, 17, 19, 23\}) 
\neq 0$.  Thus, it follows that --- assuming the Elliott-Halberstam conjecture --- 
there are infinitely many gaps 
between primes that are at most $16$.

\vskip .1 in

What can we recover unconditionally?  We are so close to proving Theorem 2 
unconditionally, that clearly some tweaking of the argument must give Theorem 1!   
The idea here is to average over sets ${\Cal H}$.  For clarity, let us now denote 
$a(n)$ above by $a(n;{\Cal H})$ to exhibit the dependence on ${\Cal H}$. 
 
Given $\epsilon >0$ we wish to find primes $p$ between $x$ and $2x$ such that 
$p_{\text{next}} -p \le \epsilon \log x$.  This would prove Theorem 1.  Set $h=\epsilon \log x$, 
and let $k$ be a natural number chosen in terms of $\epsilon$, but fixed compared 
to $x$.   Consider the following two sums: 
$$
\sum\Sb 1\le h_1 < h_2 < \ldots < h_k \le h \endSb \ \ \sum_{x\le n\le 2x} a(n;\{h_1,\ldots,h_k\}), \tag{13}
$$
and 
 $$
 \sum\Sb 1\le h_1 < h_2 < \ldots < h_k \le h \endSb \ \  \sum_{1\le \ell \le h} 
 \sum\Sb x\le n\le 2x \\ n+\ell \text{ prime} \endSb a(n;\{h_1,\ldots,h_k\}). \tag{14}
 $$
If we could prove that (14) is larger than (13), it would follow that for some $n$ 
between $x$ and $2x$, there are two prime numbers between $n+1$ and $n+h$, as desired.

Our analysis above already gives us the asymptotics for (13) and (14).  
Using (9) we see that the quantity (13) is 
$$
\sim \frac{x}{(\log R)^k} \Big( \int_0^1 \frac{y^{k-1}}{(k-1)!} P^{(k)}(1-y)^2  dy 
\Big) \sum\Sb 1\le h_1 <h_2 <\ldots< h_k \le h\endSb {\frak S}(\{h_1,\ldots,h_k\}),
$$
and using Gallagher's result (3) this is 
$$
\sim \frac{x}{(\log R)^k} \frac{h^k}{k!}  \int_0^1 \frac{y^{k-1}}{(k-1)!} P^{(k)}(1-y)^2  dy. \tag{15}
$$

Now let us consider (14).  Here we distinguish two cases: the case when $\ell =h_j$ for 
some $j$, and the case when $\ell \neq h_j$ for all $j$.  The former case is handled by our analysis leading 
up to (11).  Upon using (3) again, these terms contribute 
$$
\sim k \frac{x}{(\log x) (\log R)^{k-1}} \frac{h^k}{k!} \int_0^1 \frac{y^{k-2}}{(k-2)!} 
P^{(k-1)}(1-y)^2 dy. 
\tag{16}
$$
If we choose $P(y)=y^{k+r}$ as before, we see that (16) is already just a shade below (15),   
so we need the slightest bit of extra help from the terms $\ell \neq h_j$ for any $j$.  
If $n+\ell$ is prime note that 
$$
\align
a(n;\{h_1,\ldots,h_k\}) &= \Big(\sum_{d|(n+h_1)\cdots(n+h_k)} \lam_d\Big)^2 \\
&= \Big( \sum_{d|(n+h_1)\cdots(n+h_k) (n+\ell)} \lam_d \Big)^2 = a(n;\{h_1,\ldots,h_k,\ell\}),\\
\endalign
$$ 
since the divisors counted in the latter sum but not the former are all larger than $n+\ell >x >R$ 
and so $\lam_d=0$ for such divisors.  This allows us to finesse the calculation by 
simply appealing to (11) again, with $k$ replaced by $k+1$ and $\{h_1,\ldots, h_k\}$ 
by $\{h_1,\ldots,h_k,\ell\}$.  Thus the latter class of integers $\ell$ contributes
$$
\sim \sum\Sb 1\le h_1 < h_2 <\ldots <h_k \le h\endSb \sum\Sb \ell =1 \\ \ell 
\neq h_j \endSb^{h} \frac{x}{(\log x)(\log R)^k} {\frak S}(\{h_1,\ldots,h_k,\ell\}) 
\int_0^1 \frac{y^{k-1}}{(k-1)!} P^{(k)}(1-y)^2 dy.
$$ 
Appealing to (3) again --- we are now summing over $k+1$ element sets but 
each set is counted $k+1$ times --- this is 
$$
\sim \frac{x}{(\log R)^k} \frac{h^k}{k!} \frac{h}{\log x} \int_0^1 \frac{y^{k-1}}{(k-1)!} 
P^{(k)}(1-y)^2 dy. 
\tag{17}
$$
This accounts for a factor of $\epsilon$ times the quantity in (15), 
and now the combined contribution of (16) and (17) 
may be made larger than (15), proving Theorem 1!

\vskip .1 in 
 
 \noindent{\bf Primes in arithmetic progressions.}  It remains to explain what is 
 meant by the Bombieri-Vinogradov theorem and the Elliott-Halberstam conjecture.  Recall that 
 we required knowledge of these estimates for primes in progressions while 
 discussing the error terms that arise while evaluating the numerator of (8).
 
Let us write
$$
\pi(x) = \text{li}(x) + E(x),
$$
where $\text{li}(x)$ stands for the `logarithmic integral' $\int_2^x \frac{dt}{\log t}$, which is 
the expected main term, and $E(x)$ stands for an `error term'.  The main term $\text{li}(x)$ is, by integration 
by parts, roughly $x/\log x$.  As for the error term $E(x)$, the standard proofs 
of the Prime Number Theorem give that for any number $A>0$ there 
exists a constant $C(A)$ such that 
$$
|E(x)| \le C(A) \frac{x}{(\log x)^A}. 
$$
The argument generalizes readily for 
primes in progressions. Given an arithmetic progression $a\pmod q$ with 
$(a,q)=1$ let us write 
$$
\pi(x;q,a) = \frac{1}{\phi(q)} \text{li}(x) + E(x;q,a),
$$
where $\text{li}(x)/\phi(q)$ is the 
expected main term --- the primes are equally divided among the available residue classes --- 
and $E(x;q,a)$ is an `error term' which we would like to be small.  As with 
the Prime Number Theorem, for every $A>0$ there exists a constant $C(q,A)$ 
such that 
$$
|E(x;q,a)| \le C(q,A) \frac{x}{(\log x)^A}.
$$
We emphasize that the constant $C(q,A)$ may depend on $q$.   Therefore, 
this result is meaningful only if we think of $q$ as being fixed and let $x$ tend to $\infty$.  
In applications such a result is not very useful, because we may require $q$ not to be 
fixed, but to grow with $x$.  For example, in our discussions above 
we want to deal with primes in progressions $\pmod {[d_1,d_2]}$ 
which can be as large as $R^2$, and we'd like this to be of size $x^{\frac 12}$ 
and would love to have it be even larger.    Thus the key issue while discussing 
primes in arithmetic progressions is the uniformity in $q$ with 
which the asymptotic formula holds.   

What is known about $\pi(x;q,a)$ for an individual modulus $q$ is 
disturbingly weak.  From a result of Siegel we know that for any 
given positive numbers $N$ and $A$, there exists a constant $c(N,A)$ 
such that if $q < (\log x)^N$ then 
$$
|E(x;q,a)| \le c(N,A) \frac{x}{(\log x)^A}. 
$$
This is better than the result for fixed $q$ mentioned 
earlier, but the range of $q$ is still very restrictive.  An 
additional defect is that the constant $c(N,A)$ 
cannot be computed explicitly\footnote{This is not due to laziness, but is 
a fundamental defect of the method of proof.} in terms of $N$ and $A$.    

If we assume the 
Generalized Riemann Hypothesis (GRH) then we would fare much better: if $x\ge q$ there 
exists a positive constant $C$ independent of $q$ such that  
$$
|E(x;q,a)| \le C x^{\frac 12} \log x.
$$ 
This gives a good asymptotic formula for $\pi(x;q,a)$ in the range $q\le x^{\frac 12}/(\log x)^3$, say.

Given a modulus $q$ let us define 
$$
E(x;q) = \max_{(a,q)=1} |E(x;q,a)|.
$$
We have discussed above the available weak bounds for $E(x;q)$, 
and the unavailable strong GRH bound.  
Luckily, in many applications including ours, 
we don't need a bound for $E(x;q)$ for each individual 
$q$, but only a bound holding in an average sense as $q$ varies.  In the application 
to small gaps, we want primes 
in progressions $\pmod {[d_1,d_2]}$, but recall that we also have a sum over $d_1$, 
$d_2$ going up to $R$.   An extremely powerful result of Bombieri and Vinogradov 
gives such an average estimate for $E(x;q)$.  Moreover, this average result is nearly as good as 
what would be implied by the GRH.  

\proclaim{The Bombieri-Vinogradov theorem} For any positive constant $A$ 
there exist constants $B$ and $C$ such that 
$$
\sum_{q\le {\Cal Q}} \max_{y\le x} |E(y;q)| \le  C \frac{x}{(\log x)^A}, \tag{18}
$$ 
with ${\Cal Q}=x^{\frac 12}/(\log x)^B$.  
\endproclaim

The constant $B$ can be computed explicitly; for 
example $B=24A+46$ is permissible, but the constant $C$ here cannot be 
computed explicitly (a defect arising from Siegel's theorem mentioned above).    
The Bombieri-Vinogradov theorem tells us that 
on average over $q\le {\Cal Q}$ we have $E(x;q)\le C x(\log x)^{-A}/{\Cal Q} 
=C x^{\frac 12} (\log x)^{B-A}$.  Apart from the power of $\log x$, this 
is as good as the GRH bound!  

A straight-forward application of the Bombieri-Vinogradov theorem 
shows that as long as $R^2 \le x^{\frac 12}/(\log x)^B$ for 
suitably large $B$, the error terms arising in the Goldston-Pintz-\Yildirim 
argument will be manageable.  If we wish to take $R$ larger, then 
we must extend the range of ${\Cal Q}$ in (18).  Such extensions 
are conjectured to hold, but unconditionally the range in (18) has 
never been improved upon\footnote{Although, Bombieri, Friedlander and Iwaniec [4] 
have made important progress in related problems}.  

\proclaim{The Elliott-Halberstam conjecture}  Given $\epsilon >0$ 
and $A>0$ there exists a constant $C$ such that 
$$
\sum_{q\le {\Cal Q}} \max_{y\le x} |E(x;q)| \le C \frac{x}{(\log x)^A},
$$
with ${\Cal Q} =x^{1-\epsilon}$.
\endproclaim

The Elliott-Halberstam conjecture would allow us to take $R=x^{\frac 12-\epsilon}$ 
in the Goldston-Pintz-\Yildirim argument.  It is worth emphasizing 
that knowing (18) for ${\Cal Q}=x^{\theta}$ with any $\theta >\frac 12$ 
would lead to the existence of bounded gaps between large primes.  

Finally, let us mention a conjecture of Montgomery which lies deeper 
than the GRH and also implies the Elliott-Halberstam conjecture.

\proclaim{Montgomery's conjecture} For any $\epsilon>0$ there 
exists a constant $C(\epsilon)$ such that for all $q\le x$ 
we have 
$$
E(x;q) \le C(\epsilon) x^{\frac 12+\epsilon} q^{-\frac 12}.
$$
\endproclaim 

We have given a very rapid account of prime number theory.    
For more detailed accounts we refer 
the reader to the books of Bombieri [2], Davenport [6], and 
Montgomery and Vaughan [24].

 \vskip .1 in

 \noindent{\bf Future directions.} We conclude the article by mentioning a few 
 questions related to the work of Goldston-Pintz-\Yildirim\!\!.  
 
 First and 
 most importantly, is it possible to prove unconditionally the existence of bounded 
 gaps between primes?  As it stands, the answer appears to be 
 no, but perhaps suitable variants of the method will succeed.  There are 
 other sieve methods available beside Selberg's.  Does modifying 
 one of these (e.g. the combinatorial sieve) lead to a better result?   If 
 instead of primes we consider numbers with exactly two prime factors, 
 then Goldston, Graham, Pintz, and \Yildirim [13] have shown that there 
 are infinitely many bounded gaps between such numbers.  
 
 In a related  vein, assuming the Elliott-Halberstam conjecture, can one get to twin 
 primes?  Recall that under that assumption, we 
 could show that infinitely many permissible $6$-tuples contain two 
 primes.  Can the $6$ here be reduced? Hopefully, to $2$?  Again the 
 method in its present form cannot be pushed to yield twin primes, but 
 maybe only one or two new ideas are needed. 
 
Given any $\epsilon >0$, Theorem 1 shows that for 
infinitely many $n$ the interval $[n,n+\epsilon \log n]$ contains 
at least two primes.  Can we show that such intervals sometimes 
contain three primes?  Assuming the Elliott-Halberstam conjecture 
one can get three primes in such intervals, see [12].  Can this be made 
unconditional?  What about $k$ primes in such intervals 
for larger $k$?

Is there a version of this method which can be adapted to give 
long gaps between primes?  That is, can one attack Erd{\H o}s's 
 \$10,000 question?

\vskip .1 in 

{\bf Acknowledgments.}  I am very grateful to Carine Apparicio, Bryden Cais, 
Brian Conrad, Sergey Fomin, 
Andrew Granville, Leo Goldmakher, Rizwan Khan, Jeff Lagarias,  Youness Lamzouri, 
J{\'a}nos Pintz, 
and Trevor Wooley for their careful reading of this article, and many valuable 
comments.

 \enddocument